\documentclass[english]{IEEEtran}

\usepackage[T1]{fontenc}

\usepackage[letterpaper]{geometry}
\onecolumn

\ifCLASSINFOpdf
   \usepackage[pdftex]{graphicx}
   
   \graphicspath{{Plots/}} 
\else
\fi

\usepackage{amsmath}
\usepackage{amssymb}

\usepackage[usenames]{color}

\begin{document}

\title{Adaptive Finite-time and Fixed-time Control Design using Output Stability Conditions*}

\author{Konstantin~Zimenko, Denis Efimov, and~Andrey~Polyakov%

\thanks{*The material in this paper was partially presented at the 2019
IEEE Conference on Decision and Control \cite{ZimenkoCDC2019}. 
}

\thanks{Konstantin Zimenko, Denis Efimov and Andrey Polyakov are with Faculty of Control Systems and Robotics, ITMO University, 49 Kronverkskiy av., 197101 Saint Petersburg, Russia. {\tt\small (e-mail: kostyazimenko@gmail.com).}}
\thanks{Denis Efimov and Andrey Polyakov are with  Inria, Univ. Lille, CNRS, UMR 9189 - CRIStAL, F-59000 Lille, France.
        {\tt\small (e-mail: andrey.polyakov@ inria.fr, denis.efimov@inria.fr).} }}

\maketitle

\begin{abstract}
The present paper provides a sufficient condition to ensure output finite-time and fixed-time stability. Comparing with analogous researches the proposed result is less restrictive and obtained for a wider class of systems. The presented output stability condition is used for adaptive control design, where the state vector of a plant is extended by adjustable control parameters.
\end{abstract}


\IEEEpeerreviewmaketitle

\section{Introduction}

\IEEEPARstart{F}{requently}, the control practice needs regulation algorithms, which ensure output (in particular, a part of states) convergence in a finite time (i.e., the output $Y(t, x_0) = 0$ for all $t\ge T(x_0)$ and some $0\le T(x_0)< +\infty$ dependent on the initial conditions $x(0)=x_0$) or in a fixed time (i.e., $T(x_0)\leq T_{\max}$ for all initial conditions $x_0$). Such problem statements usually appear in mechanical and robotic systems, aerospace applications, particle collision systems (see, for example, \cite{Vorotnikov1998, RumyantsevOziraner1987, FradkovPogromsky1998, Fradkov1999, Hai-Ping1995, ShiriaevFradkov2000, Haddad2015}). 

The output stabilization is a rather common control issue, certain classes of identification problems and adaptive control systems may be considered in the context of output stability. For example \cite{Narendra1989}, in the case of adaptive control design a closed-loop system has a state vector extended with adjustable control parameters: $\left[x^T\; \omega^T \right]^T$, where $x\in\mathbb{R}^n$ is the state vector of the plant and $\omega\in\mathbb{R}^r$ is the adjustable control parameters vector. Thus, a standard control goal is to guarantee output (partial) stability: the states of the plant should be stabilized at the origin asymptotically or in a finite time, while the  adjustable parameters may remain just bounded. 
A similar output stabilization problem arises in state observer design, where the common dynamics of the system includes the plant and the observer states, $x$ and $\hat{x}$, respectively, while it is necessary to ensure the convergence of the state estimation error $e=x-\hat{x}$ (in nonlinear case the dynamics of $e$ may be dependent on $x$) \cite{ArcakKokotovic2001, AndrieuTarbouriech2019}.

There are a number of results devoted to output finite-time stability (OFTS) analysis. Most of them are about partial stability analysis that is a particular case of output stability (see, for example, \cite{Haddad2015, Jammazi2012, Jammazi2014}, \cite{SunShao2019}). In \cite{ZimenkoAutomatica2020} necessary and sufficient conditions for output finite-time stability are given using Lyapunov functions. However, most of these results are obtained for special classes of systems, and/or particular control problems. The present paper provides a relaxed sufficient condition to analyze OFTS or OFxTS of a wider class of models than in \cite{ZimenkoAutomatica2020}.

A finite-time (fixed-time) control (e.g., homogeneity based) can ensure useful properties such as faster convergence, higher accuracy, and better disturbance rejection (see, for example, \cite{PolyakovBook2020}, \cite{Rosier}, \cite{PolyakovAutomatica2015}). 
Adaptive finite/fixed-time control is one of the rapidly developing areas of control theory providing an appealing performance for systems with uncertainties (see, for example, \cite{Hong2006}, \cite{UtkinPoznyak2013},  \cite{GuzmanMoreno2011}, \cite{Basin2016}, \cite{SunShao2019} and references therein). 
Despite a number of available results, the basic problem of adaptive finite/fixed-time regulation still has no solution.
In this paper, based on the new OFTS and output fixed-time stability (OFxTS) conditions, a scheme of adaptive finite/fixed-time control design is presented. The proposed scheme allows to combine with an adaptive term different finite-time (fixed-time) control algorithms $u_{FTS}$ ($u_{FxTS}$) designed for systems without parametric uncertainties. The developed adaptive control scheme guarantees that for the state of uncertain systems the desired non-asymptotic convergence can be recovered. There is no requirement on persistence of excitation.

The paper is organized in the following way. Notation used in the paper is given in Section II. Section III recalls basics on output global asymptotic stability (oGAS), OFTS, OFxTS and homogeneity property. Section IV presents the main results on sufficient condition of OFTS/OFxTS and adaptive control design with numerical examples. Finally, concluding remarks are given in Section V.


\section{\label{sec:Notations} Notation}

\begin{itemize}
\item $\mathbb{R}^n$ denotes the $n$ dimensional Euclidean space with
vector norm $|\cdot|$;
\item $\mathbb{R}_+ = \{x\in\mathbb{R}: x>0\}$, where $\mathbb{R}$ is the field of real numbers,  and $\mathbb N$ is the set of natural numbers;
\item The symbol $\overline{1, m}$ is used to denote a
sequence of integers $1,...,m$;
\item A continuous function $\sigma:\mathbb{R}_+\cup\{0\}\to \mathbb{R}_+\cup\{0\}$ belongs to class $\mathcal{K}$ if it is strictly increasing and $\sigma(0) = 0$. It belongs to class $\mathcal{K}_\infty$ if it is also unbounded;
\item A continuous function $\beta: \mathbb{R}_+\cup\{0\}\times \mathbb{R}_+\cup\{0\} \to \mathbb{R}_+\cup\{0\}$ belongs to class $\mathcal{KL}$ if $\beta(\cdot,r)\in\mathcal{K}$ and $\beta(r,\cdot)$ is decreasing to zero for any fixed $r\in\mathbb{R}_+$;
\item By $DV (x)f(x)$ we denote  the derivative of the function $V$, if differentiable, in the line of the vector field $f$ and the upper Dini derivative for a locally Lipschitz continuous function $V$:
 $$
DV(x)f(x)=\lim_{t\to 0^+}\sup \frac{V[x+tf]-V(x)}{t}.
$$ 
\end{itemize}



\section{Preliminaries}
Consider a system in the form
\begin{equation}
\label{sys}
\dot x=f(x), \quad y=h(x)
\end{equation}
with states $x\in\mathbb{R}^n$ and outputs $y\in\mathbb{R}^p$. Let the system satisfy the following assumptions:

\textbf{(A.1)} The vector field $f:\mathbb{R}^n\to\mathbb{R}^n$ ensures forward existence and uniqueness of the system solutions at least locally in time, $f(0)=0$.

\textbf{(A.2)} The function $h:\mathbb{R}^n\to\mathbb{R}^p$ is continuously differentiable, $h(0)=0$.

\textbf{(A.3)} The vector field $f:\mathbb{R}^n\to\mathbb{R}^n$ is locally Lipschitz continuous on $\mathbb{R}^n \setminus\mathcal{Y}$, where $\mathcal{Y}=\{x\in\mathbb{R}^n: h(x)=0\}$.

For the initial conditions $x_0 \in \mathbb{R}^n$, let $\Phi(t, x_0)$ be a unique maximal solution of the system~(\ref{sys}) defined over an interval $[0, T_s (x_0))$ with some $T_s(x_0)\in\mathbb R_+\cup\{+\infty\}$ (the solutions are understood in the Filippov sense \cite{PolyakovBook2020}), $Y(t, x_0)=h(\Phi(t,x_0))$. 

Below we study only global stability and attractivity  properties of the system~(\ref{sys}). The local counterparts can be obtained by a direct restriction of the domain of validity for the presented conditions.
Note that the preliminaries in this subsection are based on theoretical framework of Input-to-Output Stability and uniform oGAS presented for locally Lipschitz continuous systems in \cite{Sontag1999}--\cite{SontagWang1999}. In \cite{ZimenkoAutomatica2020} the results on oGAS were extended for a wider class of dynamics, where the Lipschitz continuity may be violated on $\mathcal{Y}$ (see Assumption A.3).

\textbf{Definition 1 \cite{Dash}, \cite{Ingalls2001}} \textit{
The system~(\ref{sys}) is forward complete if for each $x_0\in\mathbb{R}^n$ it produces a solution
$\Phi(t,x_0)$ which is defined on $[0, +\infty)$, i.e., $T_s(x_0)=+\infty$}.

\textbf{Definition 2 \cite{Sontag1999}} \textit{The system~(\ref{sys}) has the unboundedness observability (UO) property if, for each  $x_0$ such that $T_s(x_0)<+\infty$, necessarily
\begin{equation}
\limsup_{t\to T_s(x_0)}\left|Y(t,x_0)\right|=+\infty.
\end{equation}
}

In other words, any unboundedness of the state vector can be observed using the output $y$. Hence, if the output
is known to be bounded (which is the case under
the output stability properties described below), then the UO property is equivalent to forward completeness \cite{Ingalls2001}. Note, that any system has UO property in the output $h(x) = x$.

\textbf{Definition 3 \cite{Sontag2000}, \cite{SontagECC}} \textit{A system~(\ref{sys}) is oGAS  if
\begin{itemize}
\item it is forward complete, and
\item there exists a $\mathcal{KL}$-function $\beta$ such that
\begin{equation}
\label{OSdef}
|Y(t, x_0)|\le \beta(|x_0|, t) \quad \forall t\ge 0
\end{equation}
holds for all $x_0 \in\mathbb{R}^n$.
\end{itemize}
If, in addition, there exists $\sigma \in\mathcal{K}$ such that
\begin{equation}
\label{OLOSdef}
|Y(t, x_0)|\le \sigma(|h(x_0)|) \quad \forall t\ge 0
\end{equation}
holds for all $x_0\in\mathbb{R}^n$, then the system is output-Lagrange output globally asymptotically stable (OLoGAS). Finally, if one
strengthens~(\ref{OSdef}) to
\begin{equation}
\label{SIOSdef}
|Y(t, x_0)|\le \beta(|h(x_0)|, t), \quad \forall t\ge 0
\end{equation}
for all $x_0\in\mathbb{R}^n$, then the system is state-independent output globally asymptotically stable (SIoGAS).}

\textbf{Lemma 1 \cite{SontagWang1999}} \textit{For system~(\ref{sys}) having the UO property, the following relations are valid:
 $$
SIoGAS \Rightarrow OLoGAS \Rightarrow oGAS.
$$ 
In the general case, all inverse relations are not satisfied.}

Let us present definitions for corresponded Lyapunov functions.

 \textbf{Definition 4 \cite{Dash}, \cite{Ingalls2001}} \textit{For the system~(\ref{sys}), a smooth function $V$ and a function $\lambda : \mathbb{R}^n \to \mathbb{R}_+\cup\{0\}$ are called respectively an oGAS-Lyapunov function and an auxiliary modulus if there exist $\xi_1, \xi_2 \in \mathcal{K}_\infty$ such
that 
\begin{equation}
\label{OS1}
\xi_1(|h(x)|)\le V(x)\le \xi_2(|x|) \quad \forall x\in\mathbb{R}^n
\end{equation}
 holds and there exists $\xi_3\in\mathcal{KL}$
such that
\begin{equation}
\label{OS2}
D V(x)f(x)\le -\xi_3(V(x),\lambda(x))
\end{equation} 
for all $x\in\mathcal{X}$, where $\mathcal{X}=\{x\in\mathbb{R}^n:V(x)>0\}$, and $\lambda$ satisfies the following conditions, either
\begin{itemize}
\item[(a)] 
 $0\le\lambda(x)\le |x|$ for all $x \in \mathbb{R}^n$,
 $\lambda$ is locally Lipschitz on the set 
$\mathcal{X}$ and satisfies
\begin{equation}
\label{smod}
D\lambda(x)f(x)\le 0
\end{equation} 
for almost all $x\in \mathcal{X}$,
\end{itemize}
or
\begin{itemize}
\item[(b)] there exists some $\theta\in\mathcal{K}$ such that
\begin{equation}
\label{theta}
\lambda(\Phi(t,x_0))\le \theta(|x_0|)
\end{equation} 
for all $t\ge 0$ and $x\in\mathcal{X}$.
\end{itemize}
The function $V$ is called an OLoGAS-Lyapunov function if it is an oGAS-Lyapunov function, and in addition, inequality~(\ref{OS1}) can be strengthened to
\begin{equation}
\label{OLOS1}
\xi_1(|h(x)|)\le V(x)\le \xi_2(|h(x)|), \quad \forall x\in\mathbb{R}^n.
\end{equation}
The function $V$ is called the SIoGAS-Lyapunov function if the inequality~(\ref{OLOS1}) is satisfied and there exists  $\xi_3 \in\mathcal{K}$ such that for all $x \in\mathcal{X}$:
\begin{equation}
\label{SIOS}
D V(x)f(x)\le -\xi_3(V(x)).
\end{equation}
}
An auxiliary modulus $\lambda$ satisfying
property (a) is called a strong auxiliary modulus,
and one satisfying property (b) is a weak auxiliary
modulus \cite{Ingalls2001}. 

Note that in the case of OLoGAS- or SIoGAS-Lyapunov function we have $\mathcal{X}=\mathbb{R}^n\setminus \mathcal{Y}$.

\textbf{Remark 1}
In \cite{Sontag1999}--\cite{SontagWang1999} all given above definitions are presented in 
the sense of uniform stability with respect to inputs $u$ for the system $\dot x=f(x,u)$, $y=h(x)$.

The following theorem gives the necessary and sufficient Lyapunov characterizations of output stability for the system~(\ref{sys}). 

\textbf{Theorem 1 \cite{ZimenkoAutomatica2020}} \textit{Suppose the system~(\ref{sys}) is UO. 
\begin{itemize}
\item[(1)] The following claims are equivalent for the system:
\begin{itemize}
\item[(a)] it is OLoGAS;
\item[(b)] it admits an OLoGAS-Lyapunov function
with a weak auxiliary modulus;
\item[(c)] it admits an OLoGAS-Lyapunov function
with a strong auxiliary modulus.
\end{itemize}
\item [(2)] The following claims are equivalent for the system:
\begin{itemize}
\item[(a)] it is SIoGAS;
\item[(b)] it admits a SIoGAS-Lyapunov function.
\end{itemize}
\end{itemize}}

\subsection{Output Finite-Time Stability}

Let us present the definition of the output finite-time stability.

\textbf{Definition 5} \textit{The system~(\ref{sys}) is said to be OFTS if it is oGAS and for any $x_0\in \mathbb{R}^n$ there exists  $0\le T_0<+\infty$ such that $Y(t,x_0)= 0$ for all $t> T_0$. The function $T(x_0)=\inf\{T_0\ge 0: Y(t,x_0)=0 \; \forall t\ge T_0\}$ is called the settling-time function.}

\textbf{Definition 6} \textit{The system~(\ref{sys}) is said to be OFxTS if it is OFTS and $\sup_{x_0\in\mathbb R^n} T(x_0)<+\infty$.}

\textbf{Definition 7} The set $M$ is said to be finite-time attractive
for~(\ref{sys}) if any solution $\Phi(t,x_0)$ of~(\ref{sys}) reaches $M$ in a
finite instant of time $t = T_M(x_0)$ and remains there $\forall t\geq T_M(x_0)$. As before, $T_M :\mathbb{R}^n \to\mathbb{R}_+\cup\{0\}$ is a settling-time function. The set $M$ is fixed-time attractive if $\sup_{x_0\in\mathbb R^n} T(x_0)<+\infty$.

The paper \cite{Haddad2015} deals with partial finite-time stability that is a particular case of OFTS: 

\textbf{Theorem 2 \cite{Haddad2015}} \textit{Consider the system
\begin{equation}
\label{PFT}
\begin{array}{ll}
\dot x_1=f_1(x_1, x_2), \quad x_1(0)=x_{10},\\
\dot x_2=f_2(x_1, x_2), \quad x_2(0)=x_{20},\\
\end{array}
\end{equation}
where $x_1\in  \mathbb{R}^{n_1}$, $x_2\in \mathbb{R}^{n_2}$ are the states, $f_1: \mathbb{R}^{n_1}\times\mathbb{R}^{n_2}\to\mathbb{R}^{n_1}$ and $f_2: \mathbb{R}^{n_1}\times\mathbb{R}^{n_2}\to\mathbb{R}^{n_2}$ are such that, for every $(x_1, x_2)\in\mathbb{R}^{n_1}\times\mathbb{R}^{n_2}$, $f_1(0, x_2)=0$ and $f_1(\cdot, \cdot)$, $f_2(\cdot, \cdot)$  are jointly 
continuous in $x_1$ and $x_2$.
 If there exist a continuously differentiable function $V : \mathbb{R}^{n_1}\times\mathbb{R}^{n_2}\to\mathbb{R}$, a class $\mathcal{K}$ functions $\alpha$ and $\beta$, a continuous function $k : \mathbb{R}_+\cup\{0\}\to \mathbb{R}_+$, a real number $\mu\in(0, 1)$ such that for $(x_1,x_2)\in \mathbb{R}^{n_1}\times \mathbb{R}^{n_2}$
\begin{equation}
\label{had2}
\alpha(|x_1|)\le V(x_1,x_2)\le\beta(|x_1|),
\end{equation}
\begin{equation}
\label{PFTcon}
DV(x_1, x_2)\left[\begin{smallmatrix}
f_1(x_1,x_2)\\f_2(x_1,x_2)
\end{smallmatrix} \right]\le -k(|x_2|)V(x_1,x_2)^\mu
\end{equation}
then for $y=x_1$ the system~(\ref{PFT}) is OFTS uniformly in $x_{20}$. Moreover, there exists 
a settling-time function $T :\mathbb{R}^{n_1}\times\mathbb{R}^{n_2}\to[0,\infty)$ such that
 $$
T(x_{10},x_{20})\le q^{-1}\left(\frac{V(x_{10},x_{20})^{1-\mu}}{1-\mu} \right), \; (x_{10},x_{20})\in\mathbb{R}^{n_1}\times\mathbb{R}^{n_2},
$$
where $q:[0,\infty)\to\mathbb{R}$ is continuously differentiable and satisfies
$$
\dot q(t)=k(|x_2(t)|), \quad q(0)=0, \quad t\geq0,
$$
and $T(\cdot,\cdot)$ is jointly continuous on $\mathbb{R}^{n_1}\times\mathbb{R}^{n_2}$.}

In the paper \cite{ZimenkoAutomatica2020} necessary and sufficient Lyapunov characterizations
of output finite-time stability are presented for the class of OLoGAS and SIoGAS systems~(\ref{sys}). The following lemma on OFTS property is used in the paper.

\textbf{Lemma 2 \cite{ZimenkoAutomatica2020}} \textit{ 
Consider a forward complete system~(\ref{sys}). Let a SIoGAS-Lyapunov function satisfies the inequality
\begin{equation}
\label{Vmucond}
D V(x)f(x) \le-cV(x)^\mu
\end{equation}
for some $c\in\mathbb{R}_+$, $\mu\in(0,1)$ and all $x\in \mathbb{R}^n\setminus \mathcal{Y}$. Then the system~(\ref{sys}) is SIoGAS and OFTS.
Moreover, the settling-time function satisfies
$
T(x) \le \frac{1}{c(1-\mu)}V(x)^{1-\mu}.
$}

Similarly, the result on fixed-time attractivity of a set can be presented:

\textbf{Lemma 3} \textit{ 
Consider a forward complete system~(\ref{sys}). Let a SIoGAS-Lyapunov function satisfies the inequality~(\ref{Vmucond}) for some $c\in\mathbb{R}_+$, $\mu> 1$ and all $x\in \mathbb{R}^n\setminus \mathcal{Y}$. Then the system~(\ref{sys}) is SIoGAS and, for every $\varepsilon\in\mathbb{R}_+$, the set $B=\{x\in\mathbb{R}^n\colon V(x)<\varepsilon\}$ is fixed-time attractive with
$
T_{\max} = \frac{1}{c(\mu-1)\varepsilon^{\mu-1}}.
$}

\textbf{Proof of Lemma 3}
Since $V(x)$ is a SIoGAS-Lyapunov function we have that $V(x)>0 $ for $h(x)\neq 0$ and $V(0)=0$. Then the claim is straightforward from~(\ref{Vmucond}) with $\mu>1$. $\qquad\qquad\qquad\qquad\qquad\qquad\qquad\qquad\qquad\qquad\qquad\qquad\;\qquad\blacksquare$

 Based on Lemmas 2, 3 and \cite{PolyakovTAC2012} the following result extends the  output Lyapunov function method providing the background for OFxTS analysis (in \cite{SunShao2019} partial fixed-time stability is studied using similar arguments):

\textbf{Corollary 1} \textit{ 
Consider a forward complete system~(\ref{sys}). Let a SIoGAS-Lyapunov function satisfies the inequality
\begin{equation}
\label{VFxT}
D V(x)f(x)\le-k_1V(x)^\mu-k_2V(x)^\nu
\end{equation}
for some $k_1, k_2\in\mathbb{R}_+$, $\mu\in(0,1)$, $\nu>1$ and all $x\in \mathbb{R}^n\setminus \mathcal{Y}$.
Then the system~(\ref{sys}) is SIoGAS and OFxTS with
$$
T(x_0) \le \frac{1}{k_1(1-\mu)}+\frac{1}{k_2(\nu-1)}.
$$
}

\textbf{Proof of Corollary 1}
Due to~(\ref{VFxT}) we have
$$
\dot V(x)\le\left\{\begin{array}{ll}
-k_1V(x)^\mu \quad \text{for}\quad V(x)\le 1\\
-k_2V(x)^\nu \quad \text{for}\quad V(x)> 1\\
\end{array}\right..
$$
Hence, for any $x_0$ such that $V(x_0)>1$ the last inequality guarantees the set $\{x\in\mathbb{R}^n: V(x)\le1 \}$ will be reached in a time $t_0\leq\frac{1}{k_2(\nu-1)}$. For  $V(x_0)\le1$ by Lemma 2 we derive $y(t,x_0) = 0$ for $t\geq \frac{1}{k_1(1-\mu)}$. Therefore, the system~(\ref{sys}) is OFxTS and $y(t,x_0) = 0$ for all $t \geq \frac{1}{k_1(1-\mu)}+\frac{1}{k_2(\nu-1)}$ and $\forall x_0\in\mathbb{R}^n$. $\blacksquare$\linebreak


\subsection{Homogeneity}

Homogeneity \cite{Zubov}, \cite{PolyakovBook2020} is an intrinsic property of an object, which remains consistent with respect to some scaling. This property provides many advantages to analysis and design of nonlinear control system (including finite-time stability studies).

For $r_i \in \mathbb R_+$, $i=\overline{1, n}$, $\rho \in\mathbb{R}_+$ and $\lambda\in\mathbb{R}_+$ define vector of weights $r=\begin{bmatrix} r_1, & \ldots, & r_n \end{bmatrix}^T$, dilation matrix $D_r(\lambda)=\text{diag} \{\lambda^{r_i}\}_{i=1}^n$ and homogeneous norm
\begin{equation}
\label{norm}
\|x\|_r=\left(\sum_{i=1}^n |x_i|^{\frac{\rho}{r_i}} \right)^{\frac{1}{\rho}}.
\end{equation}

\textbf{Definition 8 \cite{Zubov}}
\textit{ 
A function $g \colon \mathbb R^n \to \mathbb R$ (vector field $f \colon \mathbb R^n \to \mathbb R^n$) is said to be $r$-homogeneous of degree $d\in\mathbb{R}$ if 
$$\begin{array}{cc}
g(D_r(\lambda)x) = \lambda^d g(x)\\
(f(D_r(\lambda)x) = \lambda^d D_r(\lambda)f(x))
\end{array}
$$  for fixed $r$, all $\lambda > 0$ and $x \in \mathbb R^n$.
}
 
Introduce the following compact set (homogeneous sphere) $\mathbb S_r=\{x\in \mathbb R^n \colon \|x\|_r=1\}$, then for any $x\in\mathbb R^n$ there is $z\in\mathbb{S}_r$ such that $x=D_r(\lambda)z$ for $\lambda=\|x\|_r$.

\textbf{Theorem 3 \cite{BacciottiRosier2005, Rosier}} 
\textit{Let $f \colon \mathbb R^n \to \mathbb R^n$ be defined on $\mathbb R^n$ and be a continuous $r$-homogeneous vector field with degree $\nu$ ($\nu<0$). If the origin of the system
$\dot x=f(x) 
$
is locally asymptotically stable then it is globally asymptotically stable (globally finite-time stable) and for any $\mu >\max\{0,-\nu\}$ there exists a continuously differentiable Lyapunov function $V$ which is $r$-homogeneous with the degree $\mu$.
As a direct consequence, the derivative $DV(x)f(x)$ is $r$-homogeneous of degree $\mu+\nu$.}

According to \cite{Zubov, Rosier} for such a $V$ there exist constants $c_1, c_2, \bar{a}, b \in\mathbb{R}_+$,
 such that
\begin{equation}
\label{c1c2}
c_1\|x\|_r^\mu\le V(x) \le c_2\|x\|_r^\mu \quad\forall x\in\mathbb R^n,
\end{equation}
\begin{equation}
\label{vhom}
\frac{\partial V(z)}{\partial z} f(z)\le-\bar a, \quad \left|\frac{\partial V(z)}{\partial z} \right|\le b \quad \forall z\in \mathbb{S}_{r}.
\end{equation}

A nonlinear system $ \dot x=f(x,u)$ is homogeneously stabilizable with degree $\nu\in \mathbb R$ if there exists a feedback $u(x)$ such that the closed-loop system is homogeneous of degree $\nu$ and globally asymptotically stable. In this case the feedback $u(x)$ is called homogenizing of degree $\nu$.



\section{Main result}

\subsection{On Sufficient Condition for Output Finite-Time Stability}

Consider the system in the form
\begin{equation}
\label{sys11}
\dot x=f(x), \quad y=h(x)
\end{equation}
where $x\in\mathbb{R}^n$ is state vector, $y\in\mathbb{R}^p$ is output, 
the vector field $f:\mathbb{R}^n\to\mathbb{R}^n$ ensures forward existence and uniqueness of the system solutions at least locally in time,
the function $h:\mathbb{R}^n\to\mathbb{R}^p$ is continuous, $f(0)=0$ and $h(0)=0$.

\textbf{Remark 2}
Note that  assumptions A.2 and A.3 are relaxed for the system~(\ref{sys11}). Thus, the system under consideration  is of a wider class than in \cite{ZimenkoAutomatica2020}--\cite{SontagWang1999}.

The following theorem provides a sufficient condition for output finite-time stability of the system~(\ref{sys11}). In contrast to all previously given results, which utilize OLoGAS and SIoGAS Lyapunov functions, here a more general oGAS scenario is considered.

\textbf{Theorem 4} \textit{Let there exist differentiable on $\mathbb{R}^n\setminus \mathcal{Y}$ functions $U:\mathbb R^n\to\mathbb R_+\cup\{0\}$ and $W:\mathbb R^n\to\mathbb R_+\cup\{0\}$ such that for $\xi_1, \xi_2\in\mathcal{K}_\infty$ the following conditions are satisfied
\begin{equation}
\label{LF1}
\xi_1(|h(x)|)\le U(x)\le \xi_2(|h(x)|), \quad \forall x\in\mathbb{R}^n,
\end{equation}
\begin{equation}
\label{LF3}
V(x)=U(x)+W(x),
\end{equation}
\begin{equation}
\label{LF4}
D V(x)f(x)\le -a U(x)^\alpha, \quad \forall x\in\mathbb{R}^n\backslash\mathcal{Y},
\end{equation}
\begin{equation}
\label{LF5}
|DW(x)f(x)|\leq \sum_{i=1}^N b_i U(x)^{\beta_i}, \quad \forall x\in\mathbb{R}^n\backslash\mathcal{Y},
\end{equation}
where $a, b_i\in\mathbb{R}_+$, $\beta_i>\alpha$, $\alpha\in(0,1)$, $i=\overline{1,N}$, $N\in\mathbb{N}$.
Then the system~(\ref{sys11}) is OFTS provided that it is UO.	
}


\textbf{Proof of Theorem 4}
Since the inequality~(\ref{LF4}) is satisfied, then $DV(x)f(x)\leq0$ that due to non-negative definiteness of $W$ and~(\ref{LF1}) implies boundedness of the output $y$. Hence, by introduced assumptions the trajectories of the system are unique and defined for all $t\geq0$.

For some sufficiently small $\vartheta\in\mathbb{R}_+$ the inequalities~(\ref{LF3})-(\ref{LF5}) imply that
\begin{equation}
\label{3211}
\begin{array}{ll}
D U(x)f(x)&\le -aU(x)^\alpha-DW(x)f(x)\\
&\le -aU(x)^\alpha+ \sum_{i=1}^N b_i U(x)^{\beta_i}\\
&<0
\end{array}
\end{equation}
for all $x\in\mathcal{A}=\{\mathbb{R}^n\backslash\mathcal{Y}:U(x)\leq \vartheta\}$ due to $\alpha<\beta_i$, $i=\overline{1,N}$, and due to~(\ref{LF1}) the set $\mathcal{A}$ is forward invariant. On the other hand, due to (\ref{LF4}) we have $D V(x)f(x)\leq -a\vartheta^\alpha$ for all $x\in\mathbb{R}^n\setminus(\mathcal{A}\cup\mathcal{Y})$, and due to~(\ref{LF3}) the set $\mathcal{A}$ will be reached in a finite time. Hence, we have that $\lim_{t\to+\infty} U(x(t))=0$. Moreover, it has been shown that for any $\varepsilon >0$ and $\delta>0$ there is $T(\varepsilon,\delta)>0$ such that $U(x(t))\leq\varepsilon$ for all $t\geq T(\varepsilon,\delta)$ provided that $|x(0)|\leq\delta$ and, consequently, the system~(\ref{sys11}) is oGAS.

Since $\lim_{t\to+\infty} U(x(t))=0$, there exists an instant of time $\tau>0$ such that 
\begin{equation}
\label{123}
-aU(x(t))^\alpha+\sum_{i=1}^N b_i U(x(t))^{\beta_i}\leq-0.5aU(x(t))^\alpha
\end{equation}
for all $t\geq\tau$ due to $\alpha<\beta_i$. Therefore, by Lemma 2
the function $U(x)$ converges to $0$ in a finite time. Then
  the system~(\ref{sys11}) is OFTS and the settling-time  is bounded as follows:
$$
T(x_0)\le \tau+ \frac{U_\tau^{1-\alpha}}{0.5a(1-\alpha)},
$$
where $U_\tau=U(x(\tau))$.
$\quad\qquad\qquad\qquad\qquad\quad\quad\quad\qquad\qquad\qquad\qquad\qquad\qquad\quad\quad\quad\qquad\qquad\qquad\qquad\qquad\qquad\quad\quad\blacksquare$
  
  Note that choosing $W(x)=0$ the conditions of Theorem 4 became similar to Theorem 2.
  
  The following theorem provides a sufficient condition for output fixed-time stability of the system~(\ref{sys11}).
 
 \textbf{Theorem 5} \textit{Let the conditions of Theorem 4 are satisfied with~(\ref{LF4}) replaced by
 \begin{equation}
 \label{LFx}
DV(x)f(x) \leq -a_1 U(x)^{\alpha_1}-a_2 U(x)^{\alpha_2}, \quad \forall x\in\mathbb{R}^n\backslash\mathcal{Y},
 \end{equation}
 and 
 \begin{equation}
 \label{LFx2}
 W(x)\leq \sigma(\rho+|h(x)|), 
 \end{equation}
where $a_1,a_2\in\mathbb{R}_+$, $\alpha_1\in(0,1)$, $\alpha_2>1$, and $\beta_i$, $i=\overline{1,N}$ in~(\ref{LF5}) satisfies $\beta_i\in(\alpha_1,\alpha_2)$, $\sigma\in\mathcal{K}_\infty$ and $\rho\in\mathbb R_+\cup\{0\}$.
Then the system~(\ref{sys11}) is OFxTS provided that it is UO.	
}

  \textbf{Proof of Theorem 5}
    Analogously to the proof of Theorem 4 one can show that $\lim_{t\to+\infty} U(x(t))=0$, and the system~(\ref{sys11}) is oGAS. Then due to~(\ref{LFx}) we have
  $$
  \begin{array}{ll}
 DU(x)f(x)&\le -a_1 U(x)^{\alpha_1}-a_2 U(x)^{\alpha_2}-DW(x)f(x)\\
  
  &\le -a_1 U(x)^{\alpha_1}-a_2 U(x)^{\alpha_2}+ \sum_{i=1}^N b_i U(x)^{\beta_i}
  \end{array}
  $$
 for all $x\in\mathbb{R}^n\backslash\mathcal{Y}$.
Since $\alpha_2>\beta_i$ then there exists a constant $U_{\tau_1}>0$  such that 
$$
-a_1 U(x)^{\alpha_1}-a_2 U(x)^{\alpha_2}+\sum_{i=1}^N b_i U(x)^{\beta_i}<-0.5 a_2 U(x)^{\alpha_2}$$ for all $x\in \mathbb{R}^{n}\backslash\{M\}$, where 
$  M=\{x\in \mathbb{R}^n: U(x)\leq U_{\tau_1}\}$. 
Therefore,
$$
DU(x)f(x) \leq -\frac{a_2}{2} U(x)^{\alpha_2} \quad  \text{for}\quad x \in \mathbb{R}^{n}\backslash\{M\}.
$$
Hence, due to Lemma 3 
the set $M$ is fixed-time attractive with the following settling-time estimate 
$$
T_M(x_0)\le \frac{1}{0.5a_2(\alpha_2-1)U_{\tau_1}^{\alpha_2-1}}.
$$
By the same arguments, there exists $U_{\tau_2}\leq U_{\tau_1}$ such that
\begin{equation}
\label{th5123}
-a_1 U(x)^{\alpha_1}-a_2 U(x)^{\alpha_2}+\sum_{i=1}^N b_i U(x)^{\beta_i}\leq-0.5a_1 U(x)^{\alpha_1}
\end{equation}
for all $x\in B=\{x\in\mathbb{R}^n: U(x)\leq U_{\tau_2}\}$.  
The function $V$ is uniformly bounded on $M \setminus B$ due to~(\ref{LFx2}) and
its time derivative is separated from zero due to~(\ref{LFx}). On 
$M \setminus B$ we have
$$
\begin{array}{cc}
U_{\tau_2}\leq V(x)\leq k_1, \quad DV(x)f(x)\leq -k_2,
\end{array}
$$ 
where $k_1=U_{\tau_1}+\sigma \left(\rho+ \xi_1^{-1}(U_{\tau_1}) \right)$, $k_2=a_1U_{\tau_2}^{\alpha_1}+a_2U_{\tau_2}^{\alpha_2}$.
Hence, the set $B$ will be reached in a finite time $t\leq\tau_2$, where 
$$\tau_2\le \frac{1}{0.5a_2(\alpha_2-1)U_{\tau_1}^{\alpha_2-1}}+\frac{k_1-U_{\tau_2}}{k_2}.$$
 Finally, due to~(\ref{th5123}), the system~(\ref{sys11}) is OFxTS and the settling-time is bounded as follows (note that $U_{\tau_1}$, $U_{\tau_2}$, $k_1$ and $k_2$ do not depend on initial conditions, their values are completely predefined by the properties of $V$, $U$, $W$ and the system dynamics):
$$
\text{\quad\quad} T(x_0) \le  \frac{1}{0.5a_2(\alpha_2-1)U_{\tau_1}^{\alpha_2-1}}\!+\!\frac{U_{\tau_2}^{1-\alpha_1}}{0.5a_1(1-\alpha_1)}\!+\!\frac{k_1-U_{\tau_2}}{k_2}.
$$ $\qquad\qquad\qquad\qquad\qquad\qquad\qquad\qquad\quad\quad\qquad\qquad\qquad\qquad\qquad\qquad\qquad\qquad\qquad\qquad\qquad\qquad\qquad\qquad\qquad\qquad\blacksquare$


In Theorem 4, the values of $b_i$ and $\beta_i$ for $i=\overline{1,N}$ may depend on $V(x(0))$. In Theorem 5, such a dependence is admitted for $\beta_i$, $i=\overline{1,N}$ since they belong to a bounded interval $(\alpha_1,\alpha_2)$, while for $b_i$, $i=\overline{1,N}$ a dependence on $V(x(0))$ or $W(x(0))$ is allowed if there exists a uniform upper bound.

 \textbf{Remark 3} 
In general, none of the functions $U(x)$, $W(x)$, $V(x)$ is an output Lyapunov function. However, with respect to Definition 4 we have:
\begin{itemize}
\item if $V(x)$ is smooth and $W(x)\le \xi_3(|h(x)|)$ for $\xi_3\in\mathcal{K}$, then $V(x)$ is SIoGAS-Lyapunov function;
\item if  $U(x)$ is smooth and $DW(x)f(x)\ge 0$, then $U(x)$ is SIoGAS-Lyapunov function and the settling-time function in Theorem 4 is bounded by
$
T(x_0)\le \frac{U_0^{1-\alpha}}{a(1-\alpha)},
$
where $U_0=U(x_0)$.
\end{itemize}


\textbf{Example 1}
Consider the system~(\ref{sys11}) with 
$$
\begin{array}{ll}
x=\begin{bmatrix}
x_1\\
x_2\\
\end{bmatrix},\quad
f(x)=\begin{bmatrix}
-\text{sign}(x_1) |x_1|^{0.5}+x_2^2x_1\\
-|x_1|^{1.5}x_2\\
\end{bmatrix}, \;
y=x_1.
\end{array}
$$
The system admits UO property, and for $U(x)=|x_1|^{1.5}$ and $W(x)=0.75x_2^2$ the conditions~(\ref{LF1})-(\ref{LF5}) are satisfied with $DV(x)f(x)\le-1.5 U(x)^{2/3}$ and $|DW(x)f(x)|\leq 1.5x_2(0)^2 U(x)\leq 2V(x(0)) U(x)$. Then the system is OFTS.


\textbf{Example 2}
Consider the system
$$
\begin{array}{ll}
\dot x_1=-\text{sign}(x_1) |x_1|^{0.5}+2x_1\sin(x_2)-\text{sign}(x_1) x_1^2,\\
\dot x_2=|x_1|^{1.5}+\sin(x_2)\sin^2(x_1),\\
y=x_1.
\end{array}
$$
The system admits UO property (the right-hand side is bounded for a bounded value of $y$). For $U(x)=|x_1|^{1.5}$ and $W(x)=3(1+\cos(x_2))$ the conditions of Theorem 5 are satisfied with $a_1=a_2=1.5$, $b_1=6$, $\alpha_1=2/3$, $\alpha_2=5/3$ and $\beta_1=1$. Then the system is OFxTS. Note also that for this system there is no Lyapunov function purely dependent on $x_1$ guaranteeing a fixed-time convergence for this coordinate as in \cite{SunShao2019}.


\subsection{Adaptive Control Design}

The presented result can be utilized for adaptive finite-time or fixed-time control design. Consider the system 
\begin{equation}
\label{sysadg}
\begin{array}{ll}
\dot x(t)=f(x(t),u(t),\theta), \quad x(0)=x_0, \quad t\ge0,\\
\end{array}
\end{equation}
where $x(t)\in\mathbb{R}^n$ is the measurable state vector, $u(t)\in\mathbb{R}^m$ is the control input, $\theta\in\mathbb{R}^{q}$ is the vector of unknown parameters and $f:\mathbb{R}^n\times \mathbb{R}^m\times\mathbb{R}^{q}\to\mathbb{R}^n$. An adaptive control for the system~(\ref{sysadg}) can be presented in the form \cite{Narendra1989}
\begin{equation}
\label{ug}
\begin{array}{ll}
u(t)=g(x(t),\omega(t)),\\
\dot{\omega}(t)=\kappa(x(t)),
\end{array}
\end{equation}
where $\omega\in\mathbb{R}^{q}$ is the vector of adjustable control parameters, $g$ and $\kappa$ are mappings defined as $g:\mathbb{R}^n\times \mathbb{R}^{q}\to\mathbb{R}^m$ and $\kappa:\mathbb{R}^n\to\mathbb{R}^{q}$ (in contrast to \cite{Hong2006, SunShao2019} we will not assume that $q=1$ and $\kappa(x)$ is a nonnegative funciton). Then the problem of adaptive control design considered in this work is to provide output finite-time stability (finite-time partial stability) of the system~(\ref{sysadg}), (\ref{ug}) with the output $y=x$ and extended state vector $\tilde{x}=\left[x^T\; \omega^T\right]^T$ (i.e., we are interested only in convergence of $x$ to the origin).

To demonstrate how the result of Theorem 4 can be utilized for adaptive finite-time control design let us consider the system~(\ref{sysadg}) in the form
\begin{equation}
\label{sysad}
\dot{x}=Ax+B\left(\phi(x)^T\theta+u \right),
\end{equation}
where $x\in\mathbb{R}^n$, $u\in\mathbb{R}$, the pair of system matrix $A\in\mathbb{R}^{n\times n}$ and control gain matrix $B\in\mathbb{R}^{n\times 1}$ is controllable and $\phi:\mathbb{R}^n\to\mathbb{R}^{q}$ is known. Then following Theorem 4 one may
obtain:


\textbf{Theorem 6} \textit{Let $u_{FTS} :\mathbb R^{n}\to \mathbb R$ be a continuous feedback control such that the system 
\begin{equation}
\label{supr}
\dot x=Ax+Bu_{FTS}(x),
\end{equation}
is finite-time stable and $r$-homogeneous of degree $\nu<0$. Let $V_{FTS} :\mathbb R^{n}\to\mathbb R_+\cup\{0\}$ be a continuously differentiable $r$-homogeneous  of degree $\mu$ ($\mu>-\nu$) Lyapunov function for~(\ref{supr}).  
Let $|\phi(x)|<c\|x\|_r^\delta$ for some $c\in\mathbb{R}_+$ and $\delta>\nu+r_{\max}$, $r_{\max}=\max_{1\leq j\leq n} r_j$.
Then the system~(\ref{sysad}) with adaptive control 
\begin{equation}
\label{uft}
\begin{array}{ll}
u(x,\omega)=u_{FTS}(x)-\phi(x)^T \omega\\
\dot \omega=\gamma \phi(x) \left(\frac{\partial V_{FTS}(x)}{\partial x}B\right)^T
\end{array}
\end{equation}
for any $\gamma\in\mathbb{R}_+$  is finite-time stable at the origin 
and the variable $\omega$ remains bounded.
}


\textbf{Proof of Theorem 6}
According to Theorem 3 $DV_{FTS}(x)\left(Ax+Bu_{FTS}(x)\right)$ is $r$-homogeneous of degree $\mu+\nu$. Then, using the homogeneity property and~(\ref{c1c2}), (\ref{vhom}) one can obtain
$$
\begin{array}{ll}
DV_{FTS}(x)\left(Ax+Bu_{FTS}(x)\right)&=

\|x\|_r^{\mu+\nu}\frac{\partial V_{FTS}(z)}{\partial z}\left(A z+Bu_{FTS}(z)\right)\\

&\le   -a V_{FTS}(x)^\alpha
\end{array}
$$
with $z\in\mathbb{S}_r$, $0<\alpha=\frac{\mu+\nu}{\mu}<1$ and some $a\in\mathbb{R}_+$. Then choosing a candidate Lyapunov function for the system~(\ref{sysad}), (\ref{uft}) with extended state vector $\tilde{x}=\left[x^T\; \omega^T\right]^T$ as
$$
V(\tilde x)=V_{FTS}(x)+0.5\gamma^{-1}\left(|\theta-\omega| \right)^2
$$
we obtain 
$$
\begin{array}{ll}
DV(\tilde x)\left(Ax+B\left(\phi(x)^T\theta+u(\tilde{x}) \right)\right)
&\le  -a V_{FTS}(x)^\alpha+\frac{\partial V_{FTS}(x)}{\partial x}B\phi(x)^T (\theta-\omega)
-\gamma^{-1}\dot{\omega}^T(\theta-\omega)\\

&=-a V_{FTS}(x)^\alpha.
\end{array}
$$
Thus, the conditions~(\ref{LF1})-(\ref{LF4}) are satisfied with $U(\tilde x)=V_{FTS}(x)$, $W(\tilde x)=0.5\gamma^{-1}\left(|\theta-\omega| \right)^2$ and according to the proof of Theorem 4 the system is oGAS with $y=x$.

Since $DV(\tilde x)\left(Ax+B\left(\phi(x)^T\theta+u(\tilde x) \right)\right)\le 0$ and the amplitudes of $x$ and $\omega$ are bounded by the corresponding functions of initial conditions,  we have
$$
\begin{array}{ll}
\left|D W(\tilde x)\left(Ax+B\left(\phi(x)^T\theta+u(\tilde x) \right)\right)\right|
&=\left| \frac{\partial V_{FTS}(x)}{\partial x}B\phi(x)^T (\theta-\omega) \right|\\
& \le 
\epsilon
 \left| \frac{\partial V_{FTS}(x)}{\partial x}B\right| \left|\phi(x) \right|\!,
\end{array}
$$
for some $\epsilon\in\mathbb{R}_+$ dependent on initial conditions. Finally, with the use of \cite{BacciottiRosier2005} and~(\ref{c1c2}), (\ref{vhom}) we obtain
$$
\begin{array}{ll}
\left|D W(\tilde x)\left(Ax+B\left(\phi(x)^T\theta+u(\tilde x) \right)\right)\right|
&\le 
c\epsilon
 \left| \frac{\partial V_{FTS}(x)}{\partial x}B\right| \|x\|_r^\delta\\
 
& \le c\epsilon
 \left| \frac{\partial V_{FTS}(z)}{\partial z}\right|\left|B\right|\left|D_r^{-1}(\|x\|_r)\right| \|x\|_r^{\delta+\mu}\\
 
&  \le c\epsilon b \left|B\right|\max\left\{ \|x\|_r^{\delta+\mu-r_{\min}}, \|x\|_r^{\delta+\mu-r_{\max}} \right\}
\\

& \le  c\epsilon b \left|B\right|\max\left\{\!c_1^{-\beta_{\max}},c_1^{-\beta_{\min}} \!\right\}
 \max\left\{\!V_{FTS}(x)^{\beta_{\max}},V_{FTS}(x)^{\beta_{\min}}\!\right\}
\\
 
& = \tilde{b} \max\{V_{FTS}(x)^{\beta_{\max}},V_{FTS}(x)^{\beta_{\min}}\}, 
\end{array}
$$
where $z\in\mathbb{S}_r$, $D_r(\cdot)$ is the dilation matrix, $ \tilde{b}  = c\epsilon b \left|B\right|\max\left\{c_1^{-\beta_{\max}},c_1^{-\beta_{\min}} \right\}$ and $\beta_{\min}=\frac{\mu-r_{\max}+\delta}{\mu}>\alpha$, $\beta_{\max} = \frac{\mu-r_{\min}+\delta}{\mu}$ with $r_{\min} = \min_{1\leq j\leq n}r_j$. Thus, all conditions of Theorem 4 are verified. $\qquad\qquad\qquad\qquad\qquad\qquad\qquad\quad\quad\qquad\qquad\blacksquare$



Due to homogeneity property of~(\ref{supr}) the presented control scheme implies some robustness abilities of~(\ref{sysad}) as, for example, it can cancel certain non-Lipschitz disturbances.

\textbf{Example 3}
To illustrate the application of the proposed adaptive
scheme, let us consider a plant defined by
$$
\begin{array}{ll}
\dot{x}_1=x_2,\\
\dot{x}_2=\theta_1\sin(x_1x_2)+\theta_2x_2^2+u,
\end{array}
$$
where $\theta=[\theta_1\; \theta_2]^T$ is the vector of unknown constant parameters and $\phi(x)=[\sin(x_1x_2)\;\; x_2^2]^T$.

In this case~(\ref{supr}) is the double integrator system. 
According to \cite{Bhat2005} we choose the finite-time control $u_{FTS}$ in the form
$$
u_{FTS}(x)=-\left\lceil x_{2}\right\rfloor ^{\alpha}-\left\lceil \chi_{\alpha}\right\rfloor^{\frac{\alpha}{2-\alpha}},
$$
where $\chi_{\alpha}=x_{1}+\frac{1}{2-\alpha}\left\lceil x_{2}\right\rfloor^{2-\alpha}$, $\alpha\in(0,1)$ and $\left\lceil x\right\rfloor ^{\beta}=|x|^{\beta}\text{sign}(x)$. This finite-time control homogenizes the double integrator system of degree $\alpha-1<0$ with the vector of weights $r=\left[
2-\alpha,  1\right]^{T}$, and the corresponding homogeneous Lyapunov function can be chosen in the form
$$
V_{FTS}(x)=\frac{2-\alpha}{3-\alpha}\left|\chi_\alpha\right|^{\frac{3-\alpha}{2-\alpha}}+sx_2\chi_\alpha+\frac{l}{3-\alpha}|x_2|^{3-\alpha},
$$
where $l$ and $s$ are positive reals.
Then for $\alpha=0.5$, $\delta=\rho=2$ and $c=1.5$ the condition $|\phi(x)|<c\|x\|_r^\delta$ is satisfied, and according to Theorem 6 the system is finite-time stable with the use of adaptive control~(\ref{uft}) 
in the form
$$
\begin{array}{rll}
u(x, \omega) &=& -\left\lceil x_{2}\right\rfloor ^{\alpha}-\left\lceil \chi_{\alpha}\right\rfloor ^{\frac{\alpha}{2-\alpha}}-\phi(x)^T\omega, \\
\dot \omega &=& \gamma \phi(x) \left(\left\lceil\chi_{\alpha}\right\rfloor^{\frac{1}{2-\alpha}}
|x_2|^{1-\alpha}
+s\chi_\alpha+(s+r)x_2|x_2|^{1-\alpha}\right).\\
\end{array}
$$
The results of simulation are
shown in Fig. 1 for $\theta=[3 \;-2]^T$, $\gamma=l=s=1$. The results of simulation with using
the logarithmic scale are shown in Fig. 2 in order to demonstrate finite-time
convergence rate of $|x|$. The transients for the control $u(x)=u_{FTS}(x)$ are shown in Fig. 3 indicating that the control without adaptive term may not guarantee stability of the system.

\begin{figure} [h]
\begin{center}
\includegraphics [width=12.2cm] {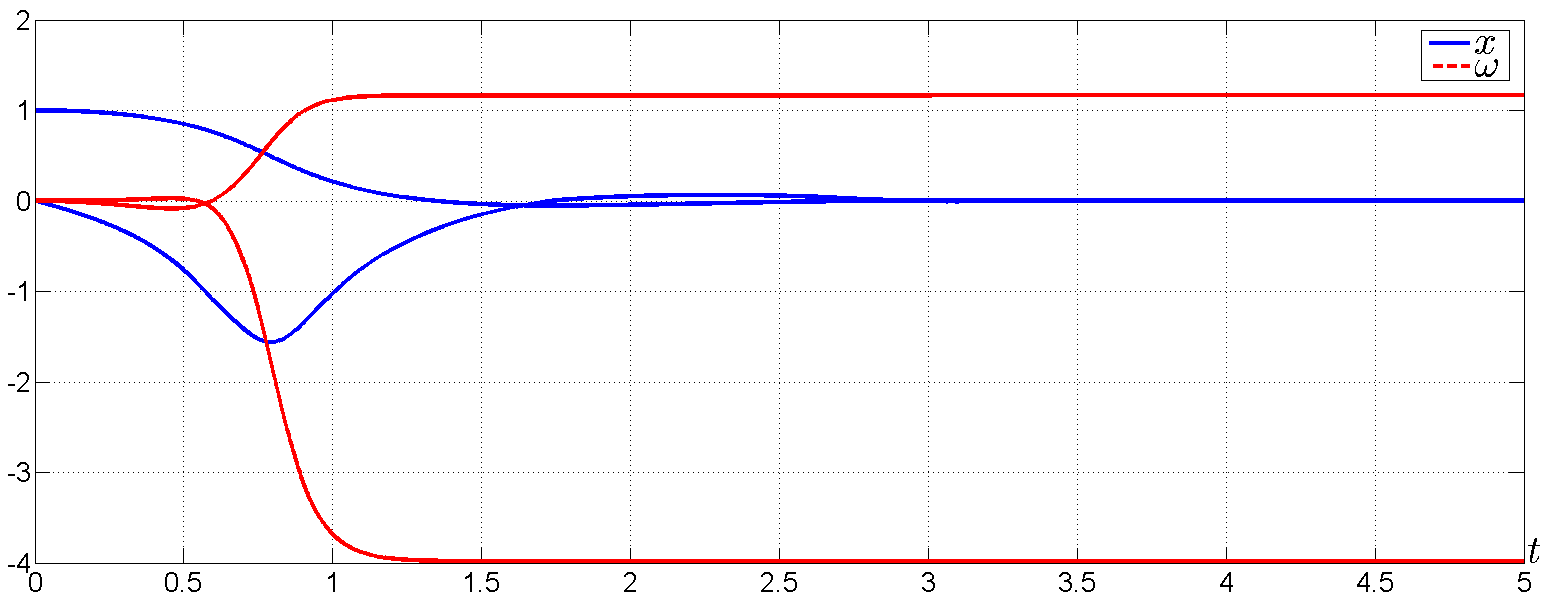}
\caption{System states $x$, $\omega$ versus time} 
\label{xom}
\end{center}
\end{figure}

\begin{figure} [h]
\begin{center}
\includegraphics [width=12.2cm] {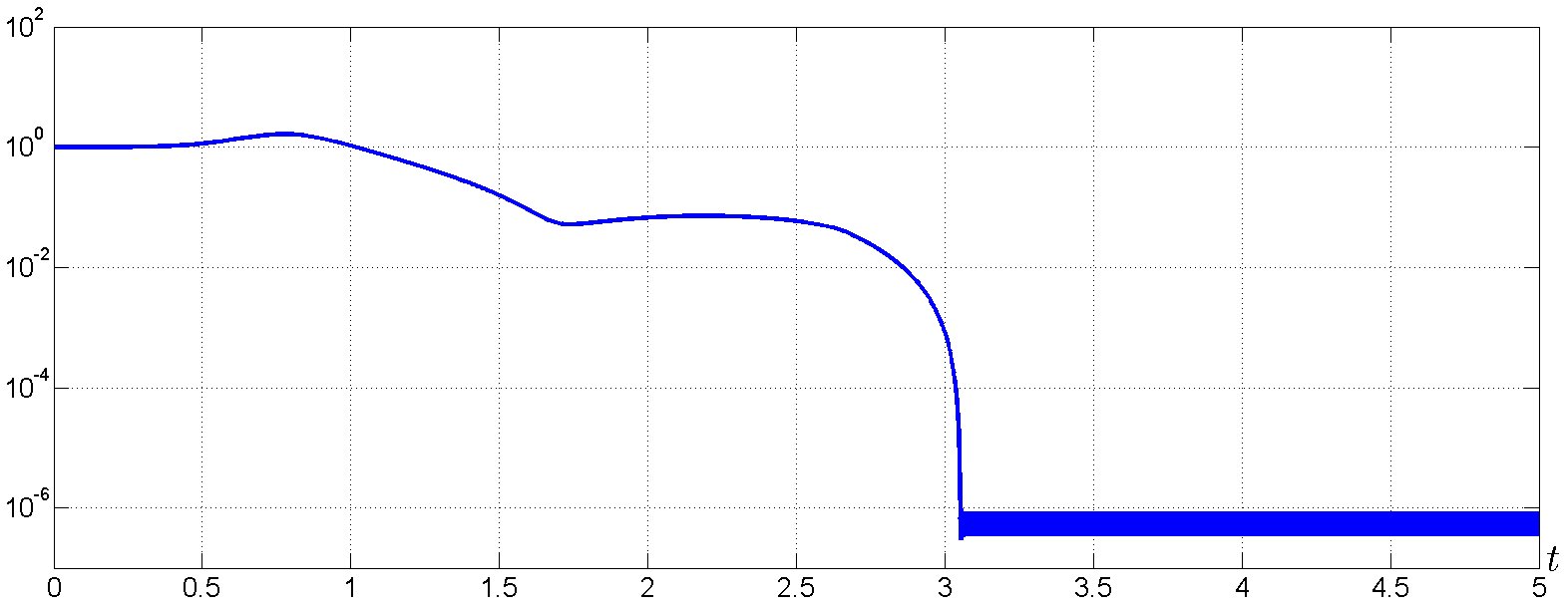}
\caption{Simulation plot of $|x|$} 
\label{xlog}
\end{center}
\end{figure}

\begin{figure} [h]
\begin{center}
\includegraphics [width=12.2cm] {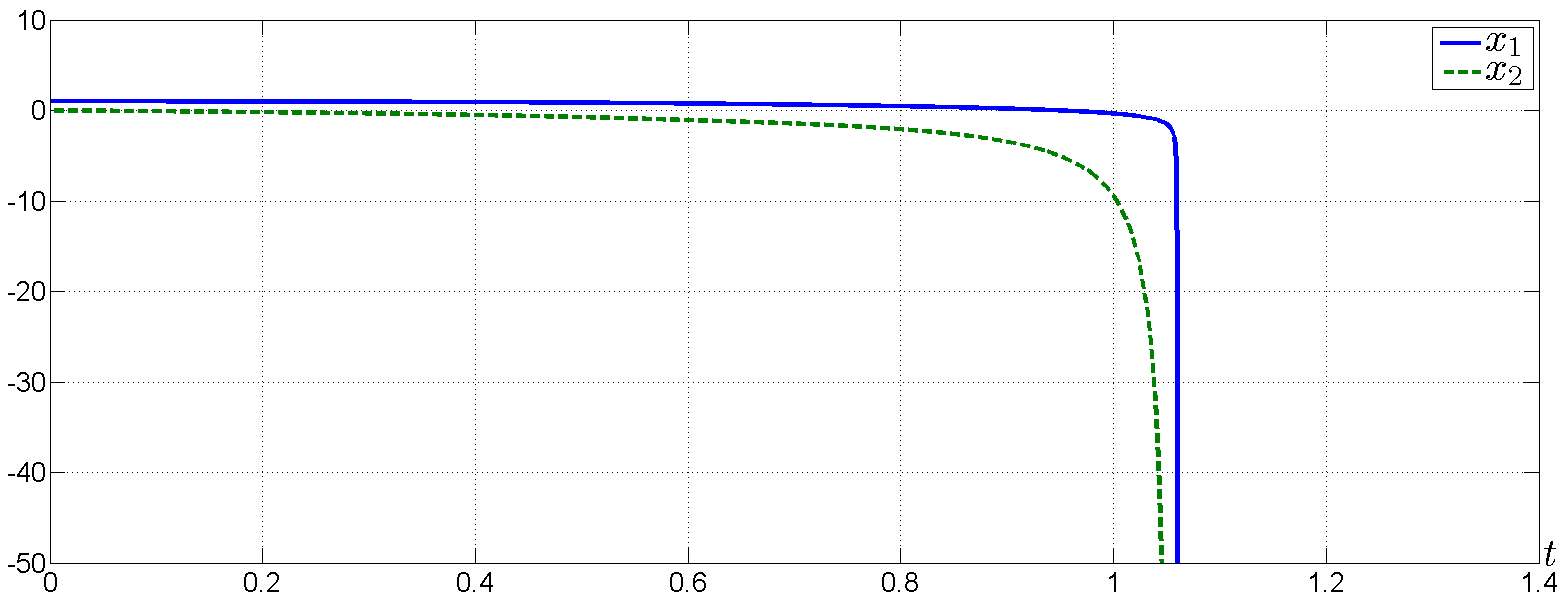}
\caption{System states $x$ for the control without adaptive term} 
\label{xust}
\end{center}
\end{figure}


\textbf{Remark 4} In \cite{PolyakovIJRNC2019} it was shown that for a stable homogeneous system $\dot{x}=f(x)$ there exists an implicitly defined homogeneous Lyapunov function $Q(V_{FTS}, x)=0$, where
$$
Q(V_{FTS}, x)=\Psi(x)^TD_r(V_{FTS}^{-1})PD_r(V_{FTS}^{-1})\Psi(x)-1,
$$
$\Psi$ is diffeomorphism on $\mathbb{R}^n\setminus\{0\}$, the homeomorphism on $\mathbb{R}^n$, $\Psi(0) = 0$. In this case by means of the implicit function theorem \cite{CourantJohn2000} we have
$$
\frac{\partial V_{FTS}}{\partial x}=-\left[\frac{\partial Q}{\partial V_{FTS}} \right]^{-1} \frac{\partial Q}{\partial x} \quad \text{for}\quad Q(V_{FTS},x)=0.
$$
The implicit Lyapunov function method is well established for homogenizing finite-time control design. For example, in order to choose homogenizing control for the system~(\ref{supr}) the results of \cite{PolyakovAutomatica2015}, \cite{PolyakovIJRNC2016} may be used.

The presented results can be extended to linear geometric homogeneous systems using their equivalence to standard homogeneous ones \cite{PolyakovIJRNC2019}.

\textbf{Remark 5} The control scheme~(\ref{uft}) presented in Theorem 6 can be applied with not necessary homogenising control laws $u_{FTS}$. Since for finite-time stable system there exists a Lyapunov function that $DV_{FTS}(x)(Ax+Bu_{FTS}(x)) \le -aV_{FTS} (x)^\alpha$, $a\in\mathbb{R}_+$, $\alpha\in(0,1)$ then the main condition for applying the scheme~(\ref{uft}) is 
$$ \left| \frac{\partial V_{FTS}(x)}{\partial x}B\right| \left|\phi(x) \right|\le \sum_{i=1}^N b_i V_{FTS}(x)^{\beta_i},$$
where $b_i\in\mathbb{R}_+$, $\beta_i>\alpha$, $i=\overline{1,N}$ and $N\in\mathbb{N}$.


Let the system~(\ref{sysad}) satisfy the following assumption:

\textbf{(A.4)} The unknown parameters are from a compact set, i.e., $|\theta| \leq \theta_{\max}$ for a known $\theta_{\max}\in\mathbb{R}_+$.

Then the results similar to Theorem 6 can be obtained for Theorem 5 by replacing a finite-time stabilizing control $u_{FTS}$ with a fixed-time one:


\textbf{Theorem 7} \textit{Let assumption A.4 be satisfied and
$u_{FxTS} :\mathbb R^n\to\mathbb R$ be stabilizing in a fixed time control law for the  system
\begin{equation}
\label{supr2}
\dot x=Ax+Bu_{FxTS}(x)
\end{equation}
and $V_{FxTS}(x)$ be the corresponding Lyapunov function satisfying 
$$
\begin{array}{lc}
D V_{FxTS}(x)(Ax+Bu_{FxTS}(x))
\leq -a_1 V_{FxTS}(x)^{\alpha_1}-a_2 V_{FxTS}(x)^{\alpha_2}, 
\end{array}
$$
$$
\xi_1(|x|)\le V_{FxTS}(x)\leq \xi_2(|x|) 
$$
with  $a_1,a_2\in\mathbb{R}_+$, $\alpha_1\in(0,1)$, $\alpha_2>1$ and $\xi_1, \xi_2\in\mathcal{K}_\infty$.
Then the system~(\ref{sysad}) with adaptive control 
\begin{equation}
\label{ufx}
\begin{array}{ll}
u(x,\omega)\!=\!u_{FxTS}(x)\!-\!\theta_{\max}\phi(x)^T\!\! \left[\arctan (\omega_1),...,\arctan (\omega_{q})\right]^T\\
\dot \omega=\gamma \theta_{\max}^{-1} \text{diag}\{1+\omega_i^2\}_{i=1}^{q}\phi(x) \left(\frac{\partial V_{FxTS}(x)}{\partial x}B\right)^T
\end{array}
\end{equation}
for any $\gamma\in\mathbb{R}_+$  is fixed-time stable at the origin
if 
\begin{equation}
\label{fxc}
\left| \frac{\partial V_{FxTS}(x)}{\partial x}B\right| \left|\phi(x) \right|\le \sum_{i=1}^N b_i V_{FxTS}(x)^{\beta_i}
\end{equation} 
is satisfied for $b_i\in\mathbb{R}_+$, $\beta_i\in(\alpha_1, \alpha_2)$, $i=\overline{1,N}$ and $N\in\mathbb{N}$.
}


\textbf{Proof of Theorem 7}
Choose a candidate Lyapunov function for the system~(\ref{sysad}), (\ref{uft}) with extended state vector $\tilde{x}=\left[x^T\; \omega^T\right]^T$ in the form
$$
\begin{array}{ll}
V(\tilde x)=V_{FxTS}(x) +0.5\gamma^{-1}\left|\tilde\theta \right|^2,
\end{array}
$$
where $\tilde{\theta}=\theta-\theta_{\max}\left[\arctan (\omega_1),...,\arctan (\omega_{q})\right]^T$.
Then we obtain 
$$
\begin{array}{lll}
DV(\tilde x) \left(Ax+B\left(\phi(x)^T\theta+u(\tilde x) \right)\right) 
 &\le&  -a_1 V_{FxTS}(x)^{\alpha_1}-a_2 V_{FxTS}(x)^{\alpha_2}
 +\frac{\partial V_{FxTS}(x)}{\partial x}B\phi(x)^T \tilde\theta\\
&&-\gamma^{-1}\theta_{\max}\dot{\omega}^T \text{diag}\left\{\frac{1}{1+\omega_i^2}\right\}_{i=1}^{q} \tilde{\theta}\\

&=& -a_1 V_{FxTS}(x)^{\alpha_1}-a_2 V_{FxTS}(x)^{\alpha_2}.
\end{array}
$$
Thus, the conditions~(\ref{LF1}), (\ref{LF3}), (\ref{LFx}) are satisfied with $U(\tilde x)=V_{FxTS}(x)$, $W(\tilde x)=0.5\gamma^{-1}\left|\tilde\theta \right|^2$ and the system is oGAS.
Due to assumption A.4 is satisfied, $W(\tilde x)\leq 0.5\gamma^{-1}\theta_{\max}^2\left(1 +0.5\sqrt{q}\pi\right)^2$ is globally bounded, i.e., the condition~(\ref{LFx2}) is satisfied.
Finally, by~(\ref{fxc}) the inequality~(\ref{LF5}) holds and all conditions of Theorem 5 are satisfied.  $\qquad\quad\qquad\qquad\qquad\quad\blacksquare$


The homogeneity property can be used for fixed-time control design. For example, the concept of homogeneity in bi-limit introduced in \cite{Andrieu} provides that an
asymptotically stable system is fixed-time stable if it is homogeneous of negative degree in $0$-limit and homogeneous of positive degree in $\infty$-limit. Based on this
the fixed-time convergence can be achieved by changing the homogeneity degree in hybrid algorithms (see, for example, \cite{PolyakovAutomatica2015}, \cite{Angulo}), \cite{ZimenkoIJC2018}).

\textbf{Corollary 2} \textit{Let $u_{FxTS}:\mathbb R^n \to 
\mathbb R$ be such that
$u_{FxTS}(x)=u_{1}(x)$ 
for $x\in \Omega$ and $u_{FxTS}(x)=u_2(x)$ for $x\in \mathbb R^n\setminus\Omega$, where 
$\Omega$ is neighborhood of the origin, 
$u_i: \mathbb R^n \to \mathbb R$ are $r_i$-homogeneous functions such that $u_1(x)=u_2(x)$ for  
  $x\in \partial \Omega$ (boundary of $\Omega$).
Let the system~(\ref{supr2}) be fixed-time stable and its continuously differentiable FxTS Lyapunov function $V_{FxTS}$ be such that 
it is $r_1$-homogeneous of degree $\mu_1$ for $x\in\Omega$ and $r_2$-homogeneous of degree $\mu_2$ for $x\in\mathbb{R}^n\setminus \Omega$.
Let $$|\phi(x)|<\left\{\begin{array}{ll}
\eta_1\|x\|_{r_1}^{\delta_1} \;\; \text{for}\; x\in\Omega,\\
\eta_2\|x\|_{r_2}^{\delta_2} \;\; \text{for}\; x\in\mathbb{R}^n\setminus \Omega\\
\end{array} \right.$$ for some $\eta_1, \eta_2\in\mathbb{R}_+$ and $\delta_1>\nu_1+r_{1 \max}$, $\delta_2<\nu_2+r_{2 \min}$, $r_{1\max}=\max_{1\leq j\leq n} r_{1j}$, $r_{2\min}=\min_{1\leq j\leq n} r_{2j}$.
Then with adaptive control in the form~(\ref{ufx}) the system~(\ref{sysad}) is fixed-time stable at the origin. 
}


\textbf{Proof of Corollary 2}
The proof is a direct consequence of Theorem 7 and the homogeneity property.
$\quad\qquad\qquad\qquad\blacksquare$


\textbf{Remark 6} As it is shown in the following example the given result can be extended for the case when $V_{FxTS}$ is continuously differentiable for $x \notin \{0\}\cup \partial\Omega$.

\textbf{Example 4}
Consider the system
$$
\begin{array}{ll}
\dot{x}_1=x_2,\\
\dot{x}_2=\theta_1\sin(x_1x_2)+\theta_2 x_2+u,
\end{array}
$$
where $\theta=[\theta_1\; \theta_2]^T$ is the vector of unknown constant parameters.
According to \cite{PolyakovAutomatica2015} choose $u_{FxTS}$ in the form
$$
u_{FxTS}=\left\{\begin{array}{ll}
V_{FxTS}^{1+\nu_1}k D_{r_1}(V_{FxTS}^{-1})x \quad \text{for} \;\; x^TX^{-1}x<1\\
V_{FxTS}^{1+2\nu_2}k D_{r_2}(V_{FxTS}^{-1})x \quad \text{for} \;\; x^TX^{-1}x\geq 1\\
\end{array} 
\right.
$$
with
\begin{itemize}
\item $\nu_1\in(-1,0)$, $\nu_2\in\mathbb{R}_+$, $r_1=\begin{bmatrix}
1-\nu_1 & 1\end{bmatrix}^T$, $r_2=\begin{bmatrix}
1 & 1+\nu_2 \end{bmatrix}^T$;
\item $k=YX^{-1}$, where $Y\in\mathbb{R}^{2\times 1}$, $X\in\mathbb{R}^{2\times 2}$ is a solution of linear matrix inequalities
$$
\begin{array}{cc}
AX+XA^T+BY+Y^TB^T+\zeta_1 X<0, \quad X>0,\\
\zeta_2 X\geq
X\text{diag}\{r_{1i}\}_{i=1}^2+\text{diag}\{r_{1i}\}_{i=1}^2 X>0,\\
\zeta_3 X\geq
 X\text{diag}\{r_{2i}\}_{i=1}^2+\text{diag}\{r_{2i}\}_{i=1}^2 X>0,\\ 
\end{array} 
$$
for some $\zeta_1, \zeta_2, \zeta_3 \in \mathbb{R}_+$;
\item $V_{FxTS}\in\mathbb{R}_+$ is defined implicitly by
$$
\left\{\begin{array}{ll}
Q_1(V_{FxTS}, x)=0  \quad \text{for} \;\; x^TX^{-1}x<1\\
Q_2(V_{FxTS}, x)=0 \quad \text{for} \;\; x^TX^{-1}x\geq 1\\
\end{array}\right.,
$$ 
where
$$
\begin{array}{ll}
Q_1(V_{FxTS},x)=x^TD_{r_1}(V_{FxTS}^{-1})X^{-1}D_{r_1}(V_{FxTS}^{-1})x-1,\\
Q_2(V_{FxTS},x)=x^TD_{r_2}(V_{FxTS}^{-1})X^{-1}D_{r_2}(V_{FxTS}^{-1})x-1.
\end{array}
$$
In order to find $V_{FxTS}$ the numerical procedures can be used (for example,
the bisection method may be utilized (see, e.g. \cite{PolyakovIJRNC2016}).
\end{itemize} 
The control $u_{FxTS}$ homogenizes the system~(\ref{supr2}) of degree $\nu_1<0$ for $x^TX^{-1}x<1$ and $\nu_2>0$ for $x^TX^{-1}x\geq 1$. 
Note that $V_{FxTS}$ is continuously differentiable for $x\notin  \{0\}\cup \{x\in\mathbb{R}^n: V(x)=1\}$. According to \cite{PolyakovAutomatica2015} the inequality
$$
\begin{array}{ll}
DV_{FxTS}(x)(Ax+Bu_{FxTS}(x))
 \leq \left\{\begin{array}{ll}
-\frac{\zeta_1}{\zeta_2} V_{FxTS}^{1+\nu_1}(x) \;\; \text{for} \; V_{FxTS}(x)<1,\\
-\frac{\zeta_1}{\zeta_3} V_{FxTS}^{1+\nu_2}(x) \;\; \text{for} \; V_{FxTS}(x)>1,\\
-\min\{\frac{\zeta_1}{\zeta_2}, \frac{\zeta_1}{\zeta_2}\} \;\; \text{for} \; V_{FxTS}(x)=1,\\
\end{array} \right.
\end{array}
$$
holds for almost all $t$ such that $x(t)\neq0$.
Then according to Theorem 7 the system is fixed-time stable with the use of adaptive control in the form~(\ref{ufx})
$$
\begin{array}{rl}
u(x,\omega)\!&=\!u_{FxTS}(x)-\theta_{\max}\phi(x)^T \left[\arctan (\omega_1)\;\arctan (\omega_2)\right]^T\\
\dot \omega&=\gamma \theta_{\max}^{-1} \text{diag}\{1+\omega_i^2\}_{i=1}^2\phi(x) \left(\frac{\partial V_{FxTS}(x)}{\partial x}B\right)^T,
\end{array}
$$
where according to Remark 4
 $$
\begin{array}{ll}
\frac{\partial V_{FxTS}(x)}{\partial x}= -2V_{FxTS}\left\{\begin{array}{ll}
\varsigma_1(x)  \quad \text{for} \;\; x^TX^{-1}x<1\\
\varsigma_2(x) \quad \text{for} \;\; x^TX^{-1}x\geq 1\\
\end{array} \right.\\
\end{array}
$$ with  
$$
\begin{array}{ll}
\varsigma_1(x)= &x^T D_{r_1}(V_{FxTS}^{-1})\left(\text{diag}\{r_{1i}\}_{i=1}^2 X^{-1}+X^{-1}\text{diag}\{r_{1i}\}_{i=1}^2\right)
 D_{r_1}(V_{FxTS}^{-1})xx^T D_{r_1}(V_{FxTS}^{-1})X^{-1} D_{r_1}(V_{FxTS}^{-1}),\\
\varsigma_2(x)= &x^T D_{r_2}(V_{FxTS}^{-1})\left(\text{diag}\{r_{2i}\}_{i=1}^2 X^{-1}+X^{-1}\text{diag}\{r_{2i}\}_{i=1}^2\right)
 D_{r_2}(V_{FxTS}^{-1})xx^T D_{r_2}(V_{FxTS}^{-1})X^{-1} D_{r_2}(V_{FxTS}^{-1}).\\
\end{array}
$$
 
The results of simulation are shown in Fig. 4 for $\theta=[3\; 2]^T$ and $x_0=[0\; 1]^T$. The results of simulation with using
the logarithmic scale are shown in Fig. 5 for different initial conditions. They show uniformity of the convergence time on the initial conditions. 
 The transients for the control $u(x)=u_{FxTS}(x)$ are shown in Fig. 6.

\begin{figure} [h]
\begin{center}
\includegraphics [width=12.2cm] {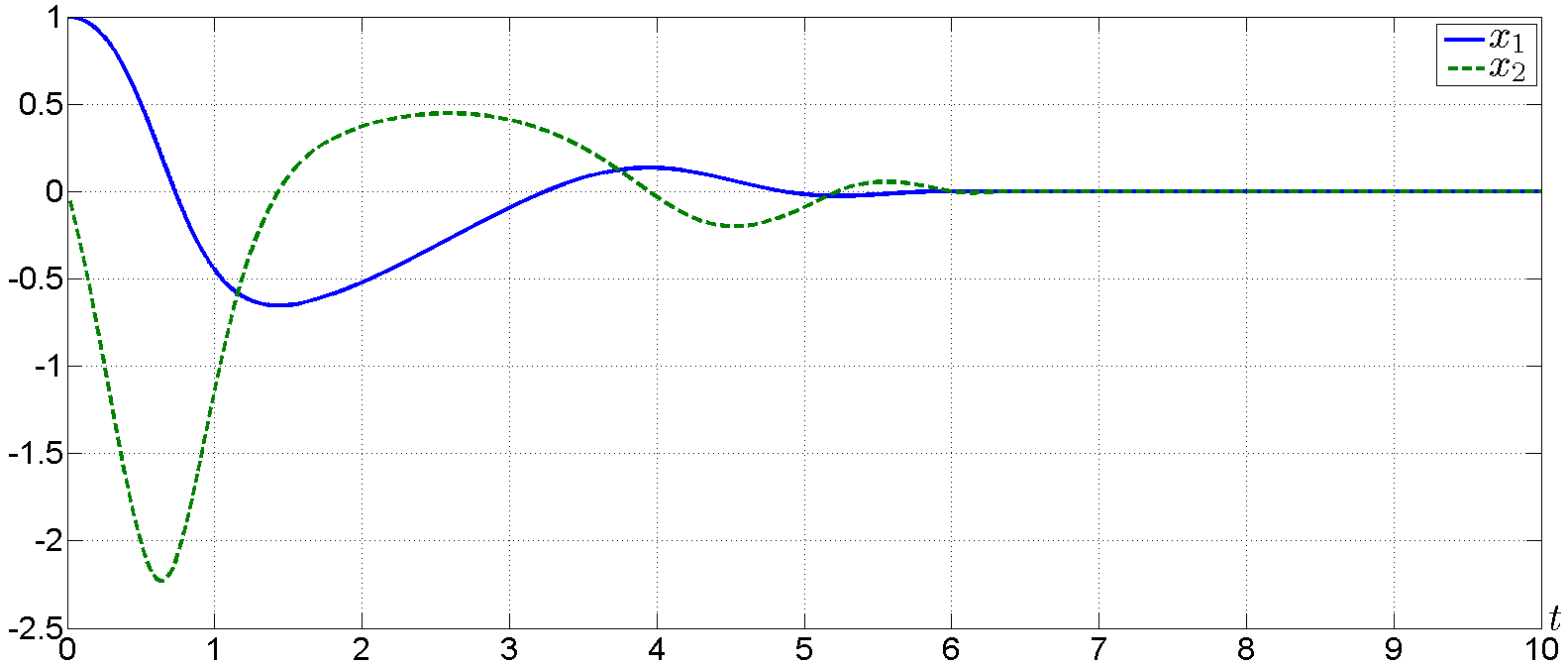}
\caption{Simulation plot for $x_0=[0\; 1]^T$} 
\label{Fx}
\end{center}
\end{figure}

\begin{figure} [h]
\begin{center}
\includegraphics [width=12.2cm] {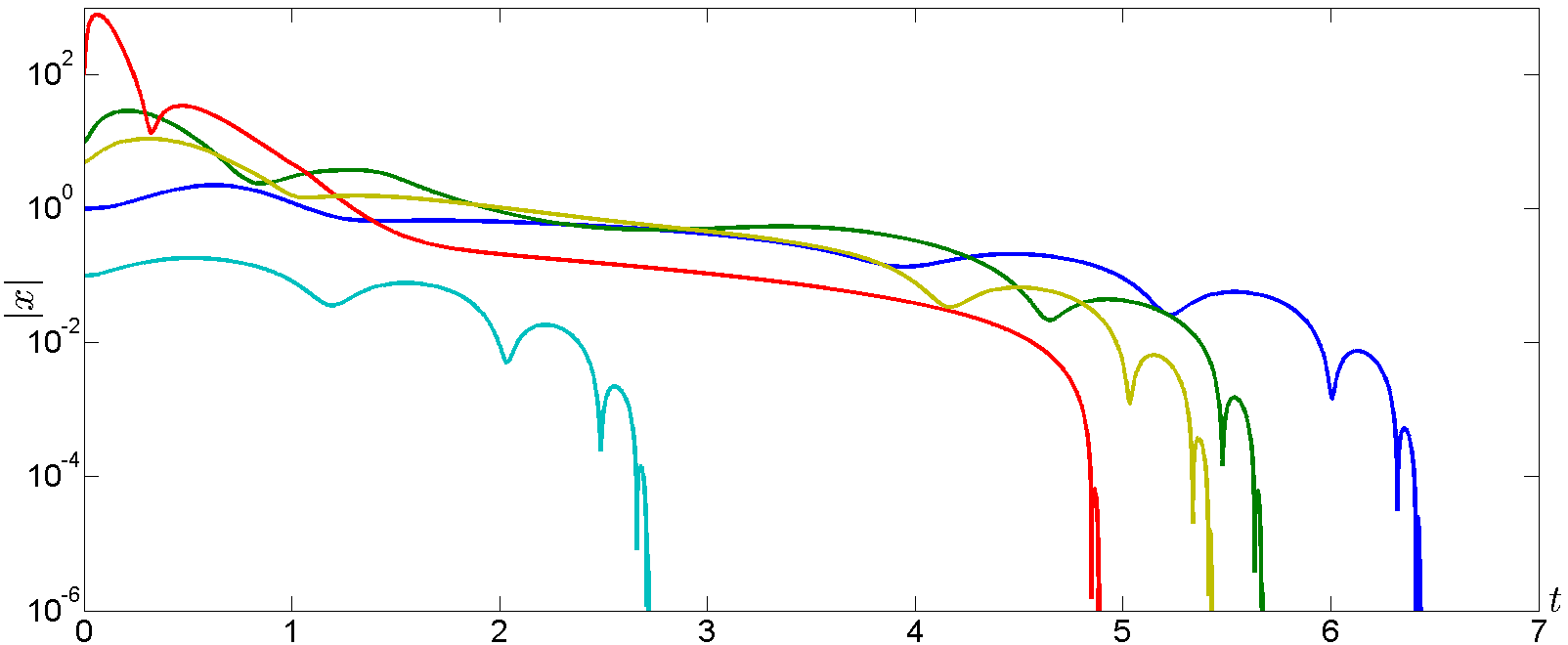}
\caption{Simulation plot for different initial conditions $x_0$} 
\label{Fxlog}
\end{center}
\end{figure}

\begin{figure} [h]
\begin{center}
\includegraphics [width=12.2cm] {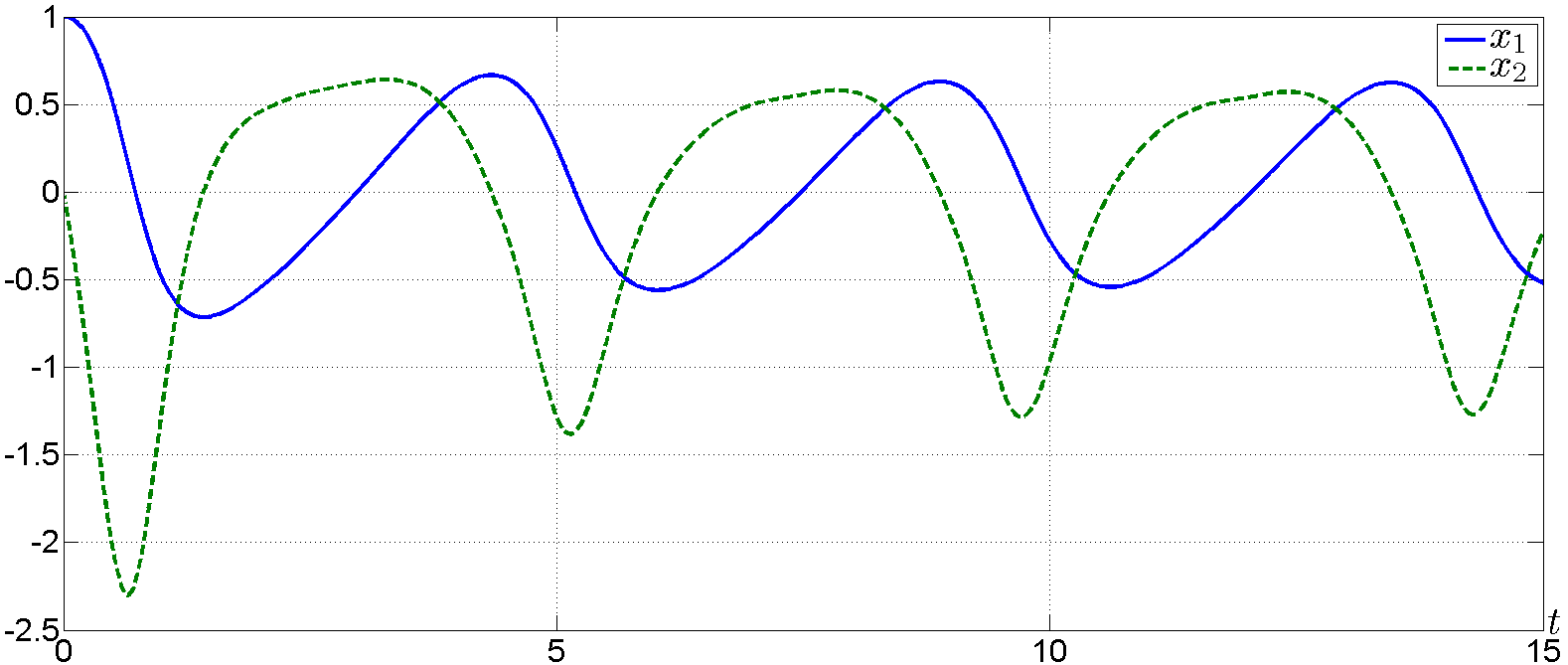}
\caption{Simulation plot for the control $u(x)=u_{FxTS}(x)$} 
\label{Fxnonad}
\end{center}
\end{figure}


\section{Conclusions}

In the paper a sufficient condition of output finite-time and fixed-time stability is presented. Comparing with existing results the presented approach is less restrictive and/or obtained for a wider class of systems. Based on the provided sufficient condition, a simple scheme of adaptive finite/fixed-time control design is presented. Possible directions for future research include control and observer design based on the use of the presented OFTS/OFxTS condition.


\begin{thebibliography}{1}


\bibitem{ZimenkoCDC2019} K.~Zimenko, D.~Efimov and A.~Polyakov, \emph{On Condition for Output Finite-Time Stability and Adaptive Finite-Time Control Scheme}, \hskip 1em plus
  0.5em minus 0.4em\relax Proc. 58th IEEE Conference on Decision and Control, Nice, 2019. 

\bibitem{Vorotnikov1998} V.I.~Vorotnikov, \emph{Partial Stability and Control}, \hskip 1em plus
  0.5em minus 0.4em\relax Birkhauser, Boston, MA,1998.

\bibitem{RumyantsevOziraner1987} V.V.~Rumyantsev and A.S.~Oziraner, \emph{Stability and stabilization of motion with respect to part of variables}, \hskip 1em plus
  0.5em minus 0.4em\relax Nauka, Moscow, 1987.

\bibitem{FradkovPogromsky1998} A.L.~Fradkov and A.Yu.~Pogromsky, \emph{Introduction to oscillations and chaos}, \hskip 1em plus
  0.5em minus 0.4em\relax World Scientific, Singapore, 1998.

\bibitem{Fradkov1999} A.L.~Fradkov, I.V.~Miroshnik and V.O.~Nikiforov, \emph{Nonlinear and adaptive control of complex systems}, \hskip 1em plus
  0.5em minus 0.4em\relax Kluwer, 1999.


\bibitem{Hai-Ping1995} Z.~Hai-Ping and M.~Feng-Xiang, \emph{On the stability of nonholonomic mechanical systems with respect to partial variables}, \hskip 1em plus
  0.5em minus 0.4em\relax Appl.Math. Mech, vol. 16(3), pp. 237--245, 1995.

\bibitem{ShiriaevFradkov2000} A.S.~Shiriaev and A.L.~Fradkov, \emph{Stabilization of invariant sets for nonlinear non-affine systems}, \hskip 1em plus
  0.5em minus 0.4em\relax Automatica, vol. 36, pp. 1709--1715, 2000.

\bibitem{Haddad2015} W.M.~Haddad and A.~L'Afflitto, \emph{Finite-time  partial stability and stabilization, and optimal
feedback control}, \hskip 1em plus
  0.5em minus 0.4em\relax Journal of the Franklin Institute, vol. 352, pp. 2329--2357, 2015.

\bibitem{Narendra1989} 
K.S.~Narendra and A.M. Annaswamy, \emph{Stable adaptive systems}, \hskip 1em plus
  0.5em minus 0.4em\relax  Prentice-Hall, Inc. Upper Saddle River, NJ, 1989.

\bibitem{ArcakKokotovic2001} M.~Arcak and P.~Kokotovi\'c, \emph{Nonlinear observers: a circle criterion design and robustness analysis}, \hskip 1em plus
  0.5em minus 0.4em\relax Automatica, vol. 37(12), pp. 1923--1930, 2001.

\bibitem{AndrieuTarbouriech2019} V.~Andrieu and S.~Tarbouriech, \emph{LMI Sufficient conditions for contraction and synchronization}, \hskip 1em plus
  0.5em minus 0.4em\relax 11th IFAC Symposium for Nonlinear Control, Vienna, 2019.


\bibitem{Jammazi2012} C.~Jammazi, \emph{A discussion on the H\"older and robust finite-time partial stabilizability of Brockett's integrator}, \hskip 1em plus
  0.5em minus 0.4em\relax ESAIM: Control, Optimisation and Calculus of Variations, vol. 18, pp. 360--382, 2012.


\bibitem{Jammazi2014} C.~Jammazi, \emph{Continuous and discontinuous homogeneous feedbacks finite-time partially stabilizing controllable multichained systems}, \hskip 1em plus
  0.5em minus 0.4em\relax SIAM Journal on Control and Optimization, vol. 52(1), pp. 520--544, 2014.

\bibitem{ZimenkoAutomatica2020} K.~Zimenko, D.~Efimov, A.~Polyakov and A.~Kremlev, \emph{On necessary and sufficient conditions for output finite-time stability}, \hskip 1em plus
  0.5em minus 0.4em\relax Automatica, vol. 125, 2021, DOI: 10.1016/j.automatica.2020.109427.


\bibitem{Sontag1999} D.~Angeli and E.D.~Sontag, \emph{Forward completeness, unboundedness observability, and their Lyapunov characterizations}, \hskip 1em plus
  0.5em minus 0.4em\relax Systems \& Control Letters, vol. 38, pp. 209--217, 1999.

\bibitem{Dash} S.N.~Dashkovskiy, D.V.~Efimov and E.D.~Sontag, \emph{Input to State Stability and Allied System Properties}, \hskip 1em plus
  0.5em minus 0.4em\relax Automation and Remote Control, vol. 72, no. 8, pp. 1579--1614, 2011.

\bibitem{Ingalls2001} B.~Ingalls and Y.~Wang, \emph{On Input-to-Output Stability for Systems not Uniformly Bounded}, \hskip 1em plus
  0.5em minus 0.4em\relax Proc.
NOLCOS, St. Petersburg, July 2001.

\bibitem{Sontag2000} E.D.~Sontag and Y.~Wang, \emph{Lyapunov Characterizations of Input to Output Stability}, \hskip 1em plus
  0.5em minus 0.4em\relax SIAM J. Control Optim., no. 39(1), pp. 226--249, 2000.

\bibitem{SontagECC} E.D.~Sontag and Y.~Wang, \emph{A Notion of Input to Output Stability}, \hskip 1em plus
  0.5em minus 0.4em\relax Proc. Eur. Control Conf., Brussels,
July 1997, DOI: 10.23919/ECC.1997.7082720.

\bibitem{SontagWang1999} E.D.~Sontag and Y.~Wang, \emph{Notions of Input to Output Stability}, \hskip 1em plus
  0.5em minus 0.4em\relax Syst. Control Lett., no. 38(4-5),
pp. 235--248, 1999.


\bibitem{Zubov} V.~Zubov, \emph{On systems of ordinary differential equations with generalized homogeneous right-hand sides (in Russian)}, \hskip 1em plus
  0.5em minus 0.4em\relax Izvestia vuzov. Mathematica, vol. 1, pp. 80--88, 1958.

\bibitem{BacciottiRosier2005} A.~Bacciotti and L.~Rosier, \emph{Lyapunov Functions and Stability in Control Theory}, \hskip 1em plus
  0.5em minus 0.4em\relax Springer, 2005.

\bibitem{Rosier} L.~Rosier, \emph{Homogeneous Lyapunov function for homogeneous continuous vector field}, \hskip 1em plus
  0.5em minus 0.4em\relax Systems \& Control Letters, vol. 19, pp. 467--473, 1992.


\bibitem{PolyakovAutomatica2015} 
A.~Polyakov, D.~Efimov and W.~Perruquetti, \emph{Finite-time and fixed-time stabilization: Implicit Lyapunov function approach}, \hskip 1em plus
  0.5em minus 0.4em\relax Automatica, vol. 51, pp. 332--340, 2015.

\bibitem{PolyakovIJRNC2016} 
A.~Polyakov, D.~Efimov and W.~Perruquetti, \emph{Robust stabilization of MIMO systems in finite/fixed time}, \hskip 1em plus
  0.5em minus 0.4em\relax International Journal of Robust and Nonlinear Control,
vol. 26(1), pp. 69-90, 2016


\bibitem{Bhat2005} 
S.~Bhat and D.~Bernstein, \emph{Continuous finite-time stabilization of the translational and rotational double integrators}, \hskip 1em plus
  0.5em minus 0.4em\relax IEEE Transactions on Automatic Control, vol. 43(5), pp. 678--682, 1998.

\bibitem{PolyakovIJRNC2019}  A.~Polyakov, \emph{Sliding Mode Control Design Using Canonical Homogeneous Norm}, \hskip 1em plus
  0.5em minus 0.4em\relax Int. J. Robust. Nonlinear Control, vol. 29(3), 682--701, 2019, DOI:10. 1002/rnc.4058.

\bibitem{CourantJohn2000} R.~Courant and F. John, \emph{Introduction to Calculus and Analysis}, \hskip 1em plus
  0.5em minus 0.4em\relax vol. II/1, New York, NY: Springer, 2000.

\bibitem{PolyakovTAC2012} A.~Polyakov, \emph{Nonlinear feedback design for fixed-time stabilization of linear control systems}, \hskip 1em plus
  0.5em minus 0.4em\relax IEEE Transactions on Automatic Control, vol. 57(8), pp. 2106-2110, 2012.

\bibitem{Andrieu}  V.~Andrieu, L.~Praly and A.~Astolfy, \emph{Homogeneous approximation, recursive observer and output feedback}, \hskip 1em plus
  0.5em minus 0.4em\relax SIAM J. Control Optim., vol. 47(4), pp. 1814-1850, 2008.


\bibitem{Angulo} M.~Angulo, J.~Moreno and L.~Fridman,  \emph{Robust exact uniformly convergent arbitrary order differentiator}, \hskip 1em plus
  0.5em minus 0.4em\relax Automatica, vol. 49, pp. 2489-2495, 2013.


\bibitem{ZimenkoIJC2018} K.~Zimenko, A.~Polyakov, D.~Efimov and W. Perruquetti, \emph{On simple scheme of finite/fixed-time control design}, \hskip 1em plus
  0.5em minus 0.4em\relax International Journal of Control, 2018, DOI: 10.1080/00207179.2018.1506889.
  
\bibitem{Hong2006} Y.G.~Hong, J.K.~Wang and D.Z.~Cheng,  \emph{Adaptive finite-time control of
nonlinear systems with parametric uncertainty}, \hskip 1em plus
  0.5em minus 0.4em\relax  IEEE Transactions on Automatic Control, vol. 51, pp. 858--862, 2006.


\bibitem{GuzmanMoreno2011}  E.~Guzman and J.A.~Moreno, \emph{A new finite-time convergent and robust direct model reference adaptive control for SISO linear
time invariant systems}, \hskip 1em plus
  0.5em minus 0.4em\relax in Decision and Control and European
Control Conference (CDC-ECC), 2011 50th IEEE Conference
on, pp. 7027–7032, 2011.



\bibitem{UtkinPoznyak2013} V.I.~Utkin, A.S.~Poznyak, \emph{Adaptive Sliding Mode Control}, \hskip 1em plus
  0.5em minus 0.4em\relax In: Bandyopadhyay B., Janardhanan S., Spurgeon S. (eds) Advances in Sliding Mode Control. Lecture Notes in Control and Information Sciences, vol. 440. Springer, Berlin, Heidelberg, 2013.
  
\bibitem{Basin2016} M.~Basin, C.B.~Panathula and Y.~Shtessel, \emph{Adaptive uniform finite-/fixed-time convergent second-order sliding-mode control}, \hskip 1em plus
  0.5em minus 0.4em\relax International Journal of Control, vol. 89(9), pp. 1777--1787, 2016, DOI: 10.1080/00207179.2016.1184759.
  
  \bibitem{SunShao2019} Z.-Y.~Sun, Y.~Shao, C.-C.~Chen, \emph{Fast finite-time stability and its application in adaptive control of high-order nonlinear system}, \hskip 1em plus
  0.5em minus 0.4em\relax Automatica, vol. 106, pp. 339--348, 2019.



\bibitem{PolyakovBook2020} A.~Polyakov, \emph{Generalized Homogeneity in Systems and Control}, \hskip 1em plus
  0.5em minus 0.4em\relax Springer, 2020.
  

 

\end{thebibliography}
\end{document}